February 26, 2020

# SUMS OF POWERS AND SPECIAL POLYNOMIALS

Khristo N. Boyadzhiev

Department of Mathematics and Statistics
Ohio Northern University
Ada, OH 45810, USA
k-boyadzhiev@onu.edu

**Abstract**

In this paper, we discuss sums of powers $1^p + 2^p + ... + n^p$ and compute both the exponential and ordinary generating functions for these sums. We express these generating functions in terms of exponential and geometric polynomials and also show their connection to other interesting series. In particular, we show their connection to an interesting problem of Ovidiu Furdui.

# 1. Introduction

Sums of powers are very popular in mathematics. Especially interesting is their connection to certain special numbers like Bernoulli and Stirling numbers – see, for example the paper of Witula et al. [8]. The present paper is also dedicated to these sums and their generating functions.

In 2014 Ovidiu Furdui published an interesting result in the form of Problem 96 in [6]. Namely, for every integer $p \geq 1$ and every real number $x$ we have

(1) $$\sum_{n=1}^{\infty} n^p \left( e^x - 1 - \frac{x}{1!} - \frac{x^2}{2!} - ... - \frac{x^n}{n!} \right) = e^x \int_0^x Q_p(t) dt,$$

where $Q_p$ is a polynomial of degree $p$, which satisfies the equation $Q_{p+1}(x) = xQ_p(x) + xQ'_p(x)$ and $Q_1(x) = x$.

Here we want to connect this result to sums of powers. It is easy to see that for every $x$,

(2) $$\sum_{n=2}^{\infty}\left(1^p+2^p+\ldots+(n-1)^p\right)\frac{x^n}{n!}=\sum_{n=1}^{\infty}n^p\left(e^x-1-\frac{x}{1!}-\frac{x^2}{2!}-\ldots-\frac{x^n}{n!}\right)$$

by changing the order of summation on the left hand side (which is easily justified)

$$\sum_{n=2}^{\infty}\left(1^p+2^p+\ldots+(n-1)^p\right)\frac{x^n}{n!}=\sum_{n=2}^{\infty}\frac{x^n}{n!}\left\{\sum_{k=1}^{n-1}k^p\right\}=\sum_{k=1}^{\infty}k^p\sum_{n=k+1}^{\infty}\frac{x^n}{n!}$$

$$=\sum_{k=1}^{\infty}k^p\left(e^x-1-\frac{x}{1!}-\frac{x^2}{2!}-\ldots-\frac{x^k}{k!}\right).$$

It is interesting to see how Problem 96 extend to the case when $p$ is not a positive integer. Equation (2) motivates us to study series as the one on the left hand side there for arbitrary $p$. We shall do this in Section 2. In the process we shall provide a new solution to Problem 96 and give more information about the polynomials $Q_p$. Later in Section 4 we shall discuss a problem similar to Problem 96 where the exponential function $e^x$ is replaced by the function $\dfrac{1}{(1-x)^m}$. This will involve the geometric polynomials defined in [3] and [5]. In the process we shall describe the exponential generating function and the ordinary generating function for the sums $1^p+2^p+\ldots+n^p$.

## 2. Sums of powers

The following two series are listed as entries 5.2.17(10) and 5.2.17 (12) in [7]

(3) $$\sum_{n=1}^{\infty}\left(1+\frac{1}{2}+\ldots+\frac{1}{n}\right)\frac{x^n}{n!}=e^x(\gamma+\ln|x|-\mathrm{Ei}(-x)),\ x\neq 0,$$

(4) $$\sum_{n=1}^{\infty}\left(1^3+2^3+\ldots+n^3\right)\frac{x^n}{n!}=e^x\frac{x}{4}(x^3+8x^2+14x+4),\ \forall x.$$

Here $\mathrm{Ei}(x)$ is the exponential integral

$$\mathrm{Ei}(x)=-\int_{-x}^{\infty}\frac{e^{-t}}{t}dt=\gamma+\ln|x|+\sum_{n=1}^{\infty}\frac{x^n}{n!n},$$

where $x\neq 0$. Clearly,

(5) $$\text{Ei}(-x) = \gamma + \ln|x| - \sum_{n=1}^{\infty} \frac{(-1)^{n-1} x^n}{n!n}, \quad \frac{d}{dx}\text{Ei}(x) = \frac{e^x}{x}, \quad \frac{d}{dx}\text{Ei}(-x) = \frac{e^{-x}}{x}.$$

We shall use also another exponential integral function

(6) $$\text{Ein}(x) = \sum_{n=1}^{\infty} \frac{(-1)^{n-1} x^n}{n!n} = \gamma + \ln|x| - \text{Ei}(-x) = \int_{x}^{\infty} \frac{1-e^{-t}}{t} dt.$$

$\text{Ein}(x)$ is an entire function as defined by the power series above which converges for every $x$.

In terms of $\text{Ein}(x)$ equation (3) becomes

(7) $$\sum_{n=1}^{\infty} \left(1 + \frac{1}{2} + \ldots + \frac{1}{n}\right) \frac{x^n}{n!} = e^x \text{Ein}(x).$$

This is the well-known exponential generating function for the harmonic numbers

$$H_n = 1 + \frac{1}{2} + \ldots + \frac{1}{n}.$$

Obviously, the two series (3) and (4) are particular cases of the general series

(8) $$M(x, p) = \sum_{n=1}^{\infty} \left(1^p + 2^p + \ldots + n^p\right) \frac{x^n}{n!},$$

where the function $M(x, p)$ is defined for all real or complex numbers $p$ and $x$. We shall see that for certain values of $p$ this series can be evaluated in closed form.

For technical convenience we define two additional entire functions

$$E(x, p) = \sum_{n=1}^{\infty} \frac{n^p x^n}{n!} \quad \text{and} \quad y_p(x) = \sum_{n=2}^{\infty} \left(1^p + 2^p + \ldots + (n-1)^p\right) \frac{x^n}{n!} = M(x, p) - E(x, p).$$

**Proposition 1**. *The following representation is true*

$$M(x, p) = E(x, p) + e^x \int_{0}^{x} e^{-t} E(t, p) dt.$$

*Proof.* Differentiating the function $y_p(x)$ with respect to $x$ we come to the differential equation

$$y_p' - y_p = E(x, p)$$

which can be written in the form

$$(e^{-x} y_p)' = e^{-x} E(x, p)$$

and therefore, its solution is

$$y_p(x) = e^x \int_0^x e^{-t} E(t, p) dt ,$$

which gives the desired representation (see also Section 7 in [4]).

We consider now several special cases. First we consider the case when $p = -1$. Then

$$e^x \int_0^x e^{-t} \left\{ \sum_{n=1}^{\infty} \frac{t^n}{n!n} \right\} dt = -e^x \int_0^x \left\{ \sum_{n=1}^{\infty} \frac{t^n}{n!n} \right\} de^{-t} = -E(x,-1) + e^x \int_0^x \frac{1-e^{-t}}{t} dt .$$

From here

$$M(x,-1) = e^x \sum_{n=1}^{\infty} \frac{(-1)^{n-1} x^n}{n!n}$$

and (3) is proved.

Next we assume that $p$ is a positive integer. In this case

(9) $$E(x, p) = \sum_{n=1}^{\infty} \frac{n^p x^n}{n!} = \left( x \frac{d}{dx} \right)^p e^x = e^x \varphi_p(x),$$

where $\varphi_p(x)$ are the exponential polynomials discussed in [1], [2], and [3]. In these papers a number of properties of $\varphi_p(x)$ were proved. The coefficients of the exponential polynomials are the Stirling numbers of the second kind $S(p,k)$, that is,

$$\varphi_p(x) = \sum_{k=0}^{p} S(p,k) x^k .$$

We have $\varphi_0(x) = 1, \varphi_1(x) = x, \varphi_2(x) = x^2 + x, \varphi_3(x) = x^3 + 3x^2 + x ,\ldots$ etc.

As pointed out in [1] and [2], these polynomials appeared as early as 1843 in the work of the prominent German mathematician Johann August Grünert (1797–1872), who founded and edited Archiv der Mathematik und Physik (for more details see [1]).

In terms of $\varphi_p(x)$ we can write

$$(10) \qquad M(x, p) = e^x \varphi_p(x) + e^x \int_0^x \varphi_p(t)dt = e^x \sum_{k=0}^p S(p,k)\left( x^k + \frac{x^{k+1}}{k+1} \right).$$

Using the property of exponential polynomials $\varphi_p(x) + \varphi'_p(x) = \varphi_{p+1}(x)/x$ (see [2], [3]) we have

$$(11) \qquad M(x, p) = e^x \int_0^x \varphi_{p+1}(t)\frac{dt}{t} = e^x \sum_{k=1}^{p+1} \frac{1}{k} S(p+1,k)x^k.$$

For $p = 3$ this is the series (4) from the beginning of the section. For $p = 1, 2$ and $p = 4$ we have correspondingly,

$$(12) \qquad \sum_{n=1}^\infty (1+2+\ldots+n)\frac{x^n}{n!} = e^x\left( \frac{x^2}{2} + x \right)$$

$$(13) \qquad \sum_{n=1}^\infty (1^2+2^2+\ldots+n^2)\frac{x^n}{n!} = e^x\left( \frac{x^3}{3} + \frac{3x^2}{2} + x \right)$$

$$(14) \qquad \sum_{n=1}^\infty (1^4+2^4+\ldots+n^4)\frac{x^n}{n!} = e^x\left( \frac{x^5}{5} + \frac{5x^4}{2} + \frac{25x^3}{3} + \frac{15x^2}{2} + x \right)$$

etc.

**Remark.** For $p = 0$ the right hand side in (10) is $e^x(1+x)$, while the left hand side is just $xe^x$. To make this formula true also for $p = 0$ we need to write the definition of $M(x, p)$ in the form

$$M(x, p) = \sum_{n=0}^\infty (0^p + 1^p + 2^p + \ldots + n^p)\frac{x^n}{n!}$$

for every integer $p \geq 0$ with the agreement $0^0 = 1$. This adjustment is not needed for the representation (11); when $p = 0$ both sides there are equal to $xe^x$.

## 3. Geometric polynomials

The exponential polynomials $\varphi_p(x)$ are related to the exponential function as demonstrated by property (9) above. We shall use also the geometric polynomials $\omega_p(x)$, $p = 0, 1, \ldots$, which are similar to the exponential polynomials and are related to the geometric series as shown in equation (16) below. They have the form

(15) $$\omega_p(x) = \sum_{k=0}^{p} S(p,k) k! x^k$$

where again $S(p,k)$ are the Stirling numbers of the second kind. Thus

$$\omega_0(x) = 1, \omega_1(x) = x, \omega_2(x) = 2x^2 + x, \ldots.$$

The geometric polynomials have the property

(16) $$\left(x \frac{d}{dx}\right)^p \frac{1}{1-x} = \sum_{n=0}^{\infty} n^p x^n = \frac{1}{1-x} \omega_p\left(\frac{x}{1-x}\right)$$

for every $p = 0, 1, \ldots$, and every $|x| < 1$ (see [1] and [3]). Note that when $p > 0$ the summation in the middle term in (16) starts, in fact, from $n = 1$. For $p = 0$ we assume that $0^0 = 1$.

Exchanging the order of summation we can immediately verify that when $p > 0$ is an integer and $|x| < 1$,

$$\sum_{n=1}^{\infty} \left(1^p + 2^p + \ldots + n^p\right) x^n = \sum_{n=1}^{\infty} x^n \left\{\sum_{k=1}^{n} k^p\right\} = \sum_{k=1}^{\infty} k^p \sum_{n=k}^{\infty} x^n$$

$$= \sum_{k=1}^{\infty} k^p x^k \sum_{n=0}^{\infty} x^n = \frac{1}{1-x} \sum_{k=1}^{\infty} k^p x^k = \frac{1}{(1-x)^2} \omega_p\left(\frac{x}{1-x}\right),$$

that is, we computed the generating function for the sums $1^p + 2^p + \ldots + n^p$, namely,

(17) $$\sum_{n=1}^{\infty} \left(1^p + 2^p + \ldots + n^p\right) x^n = \frac{1}{(1-x)^2} \omega_p\left(\frac{x}{1-x}\right), \; |x| < 1.$$

For $p = 0$ this equality is not true. To make it true for $p = 0$ we write as above

(18) $$\sum_{n=0}^{\infty} \left(0^p + 1^p + 2^p + \ldots + n^p\right) x^n = \frac{1}{(1-x)^2} \omega_p\left(\frac{x}{1-x}\right)$$

with $0^0 = 1$.

The verification of the following identities is left to the reader. Everywhere $|x| < 1$.

(19) $$\sum_{n=1}^{\infty} n^p \left(\frac{1}{1-x} - 1 - x - x^2 - \ldots - x^n\right) = \frac{x}{(1-x)^2} \omega_p\left(\frac{x}{1-x}\right)$$

for every positive integer $p$. Also,

(20) $$\sum_{n=1}^{\infty}\left(\frac{1}{1^p}+\frac{1}{2^p}+...+\frac{1}{n^p}\right)x^n = \frac{1}{1-x}\text{Li}_p(x)$$

where $\text{Li}_p(x)$ is the polylogarithm. In particular, we have the well-known series

(21) $$\sum_{n=1}^{\infty}\left(1+\frac{1}{2}+...+\frac{1}{n}\right)x^n = \frac{-\ln(1-x)}{1-x}.$$

## 4. More geometric polynomials and a result similar to Problem 96

Here we present extensions of (17) and (19). The series in (17) can be extended in the following way.

**Proposition 2.** *Let $r \geq 0$ be an integer. Then for any integer $p > 0$ and every $|x| < 1$,*

(22) $$\sum_{n=1}^{\infty}\binom{n+r}{r}\left(1^p+2^p+...+n^p\right)x^n = \frac{1}{r!}\left(\frac{d}{dx}\right)^r\left(\frac{x^r}{(1-x)^2}\omega_p\left(\frac{x}{1-x}\right)\right),$$

*where $\omega_p(x)$ are the geometric polynomials defined in Section 3.*

When $r = 0$ this is equation (17). The proof follows from the next lemma.

**Lemma 3.** *Let the function $f(x) = a_0 + a_1 x + a_2 x^2 + ...$ be analytic in a neighborhood of the origin $|x| < R$. Then for every integer $r \geq 0$ and every $|x| < R$ we have*

$$\sum_{n=0}^{\infty}\binom{n+r}{r}a_n x^n = \frac{1}{r!}\{x^r f(x)\}^{(r)}.$$

For the proof we compute directly

$$\frac{1}{r!}\{x^r f(x)\}^{(r)} = \frac{1}{r!}\left(\frac{d}{dx}\right)^r \sum_{n=0}^{\infty}a_n x^{n+r} = \frac{1}{r!}\sum_{n=0}^{\infty}a_n(n+r)(n+r-1)...(n+1)x^n$$

$$= \sum_{n=0}^{\infty}\binom{n+r}{r}a_n x^n.$$

The proposition follows from Lemma 3 and (17). We assume that $a_0 = 0^p = 0$ and take

$$f(x) = \frac{1}{(1-x)^2}\omega_p\left(\frac{x}{1-x}\right).$$

To proceed further we involve the generalized geometric polynomials $\omega_{n,r+1}(x)$ defined by

(23) $$\omega_{n,r+1}(x) = \frac{1}{r!}\sum_{k=0}^{n} S(n,k)(k+r)!\, x^k$$

for integers $r \geq 0$ and $n \geq 0$. Clearly, when $r = 0$, $\omega_{n,1}(x) = \omega_n(x)$. These polynomials have the property

(24) $$\left(x\frac{d}{dx}\right)^p \frac{1}{(1-x)^{r+1}} = \sum_{n=0}^{\infty}\binom{n+r}{r} n^p x^n = \frac{1}{(1-x)^{r+1}}\omega_{p,r+1}\left(\frac{x}{1-x}\right)$$

for any integer $p \geq 0$ and $|x| < 1$ - see [3] and [5]. Note that when $p > 0$ the summation in (24) starts from $n = 1$.

The above proposition implies the following result similar to Problem 96.

**Corollary 4.** *For every integer $r \geq 0$, every integer $p > 0$, and every $|x| < 1$ we have*

(25) $$\sum_{n=0}^{\infty} n^p \left(\frac{1}{(1-x)^{r+1}} - 1 - \binom{r+1}{1}x - \binom{r+2}{2}x^2 - \ldots - \binom{r+n}{n}x^n\right)$$

$$= \frac{1}{r!}\left(\frac{d}{dx}\right)^r \left(\frac{x^r}{(1-x)^2}\omega_p\left(\frac{x}{1-x}\right)\right) - \frac{1}{(1-x)^{r+1}}\omega_{p,r+1}\left(\frac{x}{1-x}\right).$$

When $r = 0$, equation (25) becomes equation (19). The construction of the left hand side is based on the expansion

$$\sum_{n=0}^{\infty}\binom{r+n}{n} x^n = \frac{1}{(1-x)^{r+1}}$$

for any $|x| < 1$.

*Proof.* First remember that

$$\binom{r+k}{k} = \binom{r+k}{r}.$$

Now we compute by changing the order of summation

$$\sum_{k=1}^{\infty} k^p \left(\frac{1}{(1-x)^{r+1}} - 1 - \binom{r+1}{1}x - \binom{r+2}{2}x^2 - \ldots - \binom{r+k}{k}x^k\right)$$

$$= \sum_{k=0}^{\infty} k^p \left\{ \sum_{n=k+1}^{\infty} \binom{n+r}{r} x^n \right\} = \sum_{n=1}^{\infty} \binom{n+r}{r} x^n \left\{ \sum_{k=0}^{n-1} k^p \right\}$$

$$= \sum_{n=1}^{\infty} \binom{n+r}{r} x^n \left\{ \sum_{k=1}^{n} k^p \right\} - \sum_{n=1}^{\infty} \binom{n+r}{r} n^p x^n$$

$$= \frac{1}{r!} \left( \frac{d}{dx} \right)^r \left( \frac{x^r}{(1-x)^2} \omega_p \left( \frac{x}{1-x} \right) \right) - \frac{1}{(1-x)^{r+1}} \omega_{p,r+1} \left( \frac{x}{1-x} \right)$$

and the proof is completed.

## References

[1] K. N. Boyadzhiev, *Close encounters with the Stirling numbers of the second kind,* Math. Mag., 85 (4) (2012), 252-266.

[2] K. N. Boyadzhiev, *Exponential polynomials, Stirling numbers, and evaluation of some Gamma integrals,* Abstract and Applied Analysis, Volume 2009, Article ID 168672 (electronic).
DOI:10.1155/2009/168672

[3] K. N. Boyadzhiev, *A Series transformation formula and related polynomials,* International Journal of Mathematics and Mathematical Sciences ,Vol. 2005 (2005), Issue 23, Pages 3849-3866.
DOI:10.1155/IJMMS.2005.3849

[4] K. N. Boyadzhiev, *Polyexponentials,* (2007), at  https://arxiv.org/abs/0710.1332

[5] K. N. Boyadzhiev and A. Dil, *Geometric polynomials: properties and applications to series with zeta values,* Anal. Math., 42 (3) (2016), 203–224.
DOI: 10.1007/s10476-016-0302-y

[6] O. Furdui, *Problem 96*, Mathproblems 4(2) (2014), 263. Solutions in 4(3) (2014), 304-308.

[7] A. P. Prudnikov, Yu. A. Brychkov, O. I. Marichev, Integrals and Series, vol.1 Elementary Functions, Gordon and Breach 1986.

[8]   R. Witula, K. Kaczmarek, P. Lorenc, E. Hetmaniok, and M. Pleszczyński,
*Jordan numbers, Stirling numbers, and sums of powers*, Discussiones Mathematicae,
General Algebra and Applications 34 (2014) 155–166.
DOI:10.7151/dmgaa.1225